\documentclass[preprint]{sig-alternate}

\DeclareMathOperator{\End}{End}
\DeclareMathOperator{\Inn}{Inn}
\DeclareMathOperator{\Gal}{Gal}
\DeclareMathOperator{\charak}{char}
\DeclareMathOperator{\Irr}{Irr}
\DeclareMathOperator{\id}{id}

\begin{document}
\conferenceinfo{ISSAC'07,} {July 29--August 1, 2007, Waterloo, Ontario, Canada.}
\CopyrightYear{2007}
\crdata{978-1-59593-743-8/07/0007}
\title{The Isomorphism Problem for Cyclic Algebras\\and an Application}
\numberofauthors{1}
\author{
\alignauthor Timo Hanke\titlenote{Supported by the DAAD (German Academic Exchange Service, Kennziffer D/02/00701)
and by Universidad Nacional Aut\'onoma de M\'exico, M\'exico, D.F.}\\
\affaddr{Instituto de Matemáticas, UNAM}\\
\affaddr{Circuito Exterior}\\
\affaddr{Ciudad Universitaria}\\
\affaddr{04510 México, D.F., México}\\
\email{e-mail : hanke@math.uni-potsdam.de}
}
\toappear{To appear in ISSAC'07, July 29--August 1, 2007, Waterloo, Ontario, Canada.}
\maketitle
\begin{abstract}
  The isomorphism problem means to decide if two given finite-dimensional simple algebras with center $K$ are $K$-isomorphic
  and, if so, to construct a $K$-isomorphism between them.
  Applications lie in computational aspects of 
  representation theory, algebraic geometry and Brauer group theory.

  The paper presents an algorithm for cyclic algebras that reduces the isomorphism problem to field theory
  and thus provides a solution if certain field theoretic problems including norm equations can be solved
  (this is satisfied over number fields).
  As an application, we can compute all automorphisms of any given cyclic algebra over a number field.
  A detailed example is provided which leads to the construction of an
  explicit noncrossed product division algebra.
\end{abstract}
\newtheorem{theorem}{Theorem}
\newtheorem{lemma}{Lemma}
\newtheorem{prop}{Proposition}
\newtheorem{alg}{Algorithm}
\newdef{remark}{Remark}
\newdef{rem}{Remark}
\newdef{spp}{Splitting Problem}
\newdef{ip}{Isomorphism Problem}
\newdef{ep}{Extension Problem}
\newdef{spcase}{Special Case}
\newdef{gcase}{General Case}
\newlength{\myskip}
\def\thegcase{\settowidth{\myskip}{\ }\hspace{-\myskip}\hspace{.1mm}}
\def\therem{\settowidth{\myskip}{\ }\hspace{-\myskip}\hspace{.1mm}}
\def\thespp{\settowidth{\myskip}{\ }\hspace{-\myskip}\hspace{.1mm}}
\def\theip{\settowidth{\myskip}{\ }\hspace{-\myskip}\hspace{.1mm}}
\def\theep{\settowidth{\myskip}{\ }\hspace{-\myskip}\hspace{.1mm}}
\newcommand{\comment}[1]{}
\newcommand{\set}[1]{\{#1\}}
\newcommand{\sett}[2]{\{#1\,|\,#2\}}
\newcommand{\gen}[1]{\langle #1\rangle}
\newcommand{\eps}{\varepsilon}
\newcommand{\qt}[1]{\quad\text{#1}\quad}
\newcommand{\ra}{\rightarrow}
\newcommand{\ovl}[1]{\overline{#1}}
\newcommand{\wt}{\widetilde}
\newcommand{\field}[1]{\mathbb{#1}}
\newcommand{\Q}{\field{Q}}
\newcommand{\N}{\field{N}}
\newcommand{\Z}{\field{Z}}
\newcommand{\x}{{\mathbf x}}
\newcommand{\pre}[2]{{}_{#1}#2}

\newcommand{\chieta}{
\left(
\begin{matrix}
  303\alpha^2 - 154\alpha - 276&
  314\alpha^2 + 218\alpha - 326&
  -48\alpha^2 + 151\alpha + 157
  \\
  390\alpha^2 + 708\alpha - 855&
  40\alpha^2 - 238\alpha + 430&
  -397\alpha^2 - 27\alpha + 275
  \\
  -106\alpha^2 + 25\alpha + 543&
  -128\alpha^2 - 46\alpha - 30&
  135\alpha^2 + 38\alpha - 63
  \\
\end{matrix}
\right)
} 

\newcommand{\chiw}{
\left(
\begin{matrix}
  0&
  \alpha^2 + \alpha&
  0
  \\
  0&
  -\alpha + 1&
  0
  \\
  0&
  -1&
  0
  \\
\end{matrix}
\right)
} 
\comment{
}
\section*{Introduction}

Let $K$ be a field 
and let $A_1$ and $A_2$ be two finite-dimensional central-simple $K$-algebras
($A_i$ has no proper two-sided ideal and the center is $K$).
The {\em isomorphism problem} for $A_1$ and $A_2$ means the problem to decide
whether $A_1$ and $A_2$ are $K$-isomorphic and, if so, 
to construct a $K$-isomorphism between them.
We assume that $A_1$ and $A_2$ have equal dimension,
for otherwise they are trivially non-isomorphic.

The special case when $A_2$ is the full matrix ring over $K$
is called the {\em splitting problem} for $A_1$.
We shall call a $K$-isomorphism $A_1\ra M_n(K)$ a {\em splitting} of $A_1$.

There are several applications.
For instance, to compute the irreducible representations over $K$ 
of a finite group $G$ with $|G|$ not divisible by $\charak K$,
one can decompose the semisimple group ring $KG$ into its simple components
(see \cite{eberly:decomp-nf} for an algorithm) 
and then solve the splitting problem for each component. 
As another example,
finding $K$-rational points on a Brauer-Severi variety $V$
is equivalent
to the splitting problem for the central-simple $K$-algebra
associated with $V$
(cf.\ \cite{graaf-schicho}).
The splitting problem also occurs if one wants to compute
orthogonal idempotent generators in central-simple algebras.
In this paper we pursue an application that is motivated by explicit algebra constructions. 
Because various constructions make use of automorphisms of simple algebras
that are nontrivial on the center,
we study (in section \ref{sec:ext}) 
the problem of extending an automorphism of the field $K$ 
to an automorphism of the central-simple $K$-algebra $A$,
and show that it reduces to the isomorphism problem.

The present paper solves the isomorphism problem for algebras that are presented as cyclic algebras.
A {\em cyclic algebra} is a central-simple algebra that contains a maximal subfield
(i.e.\ a subfield with maximal degree)
which is cyclic over the center.
We say an algebra is {\em presented as a cyclic algebra} if a cyclic maximal subfield is explicitly given.
For example, 
this can be a presentation by structure constants
plus an explicit generator of a cyclic maximal subfield.
Having only the theoretical information that the algebra is cyclic,
e.g.\ if the center is a global field (\cite[Thm.\ 32.20]{reiner:max-orders}),
is not sufficient.
There is no algorithm available that can produce cyclic maximal subfields of central-simple algebras
except for certain small degrees.
Indeed, for quaternion algebras this is trivial,
and for cubic algebras
one finds a cyclic maximal subfield simply by proceeding 
along the lines of Wedderburn's proof \cite[\S15.6]{pierce:ass-alg}
or Haile's proof \cite{haile:useful}
that every division algebra of degree three is cyclic.

The algorithms of this paper work by reducing the isomorphism problem to norm equations.
Norm equations are in general hard to solve,
but algorithms are known over number fields (e.g.\ Simon \cite{simon:neq}) 
and of course over finite fields.
Using computer algebra systems like KASH \cite{kash} and MAGMA \cite{magma}
for those norm equations,
our algorithms 
actually become applicable over number fields (and finite fields).

\section{Preliminaries}

Let $K$ be a field. 
Unless stated otherwise, all algebras are finite-dimensional $K$-algebras.
The tensor product $\otimes$ and the isomorphism $\cong$ without subscripts mean tensor product and isomorphism over $K$, respectively.
A $K$-algebra $A$ is called {\em central-simple}
if $A$ has no proper two-sided ideals 
and its center $Z(A)$ is $K$.
The reader is assumed to be familiar with the basic theory of central-simple algebras
as in the textbook sources Pierce \cite{pierce:ass-alg} or Reiner \cite{reiner:max-orders}.
A few relevant terms are briefly recalled in the sequel.

Let $A$ be a central-simple $K$-algebra.
The {\em degree} of $A$, denoted $\deg A$,
is the square root of the dimension $[A:K]$
(the dimension is always a square).
The algebra $A$ is called {\em split} if $A\cong M_n(K)$ where $n=\deg A$.
The {\em opposite algebra} $A^\circ$ of $A$ is the $K$-space $A$ with multiplication redefined by $a\circ b:=ba$.
If $B$ is another central-simple $K$-algebra of degree $n$ then
\begin{equation}
A\cong B\iff A\otimes B^\circ\cong M_{n^2}(K).
  \label{eq:op}
\end{equation}
A {\em maximal subfield} of $A$ is a commutative subfield $L\subseteq A$ with $[L:K]=\deg A$.
The algebra $A$ is called a {\em crossed product}
if $A$ contains a maximal subfield that is Galois over~$K$.
Moreover, $A$ is called {\em cyclic} (resp.\ {\em bicyclic}) if $A$ contains a maximal subfield
that is cyclic (resp.\ bicyclic) over $K$.
\subsection{Cyclic Algebras}\label{sec:cyclic}
Suppose $A$ is cyclic of degree $n$
and let $L$ be a maximal subfield cyclic over $K$
with $\Gal(L/K)=\gen{\sigma}$.
By the Skolem-Noether Theorem there is an element $v\in A^*$ such that
the inner automorphism $\Inn(v):x\mapsto vxv^{-1}$ of $A$ satisfies
$\Inn(v)|_L=\sigma$.
For any such $v$ we have $v^n\in K$ and
$$ A=\oplus_{i=0}^{n-1} Lv^i$$
as a $K$-space.
Setting $a=v^n$ we denote this algebra by
$$A=(L/K,\sigma,a,v)$$
(in the literature $v$ is usually omitted from the notation).
For any integer $k$ relatively prime to $n$,
\begin{equation}
(L/K,\sigma,a,v)=(L/K,\sigma^k,a^k,v^k).  
  \label{eq:k}
\end{equation}
We have 
\begin{equation*}
(L/K,\sigma,a,v)\cong M_n(K) \iff a\in N_{L/K}(L)
  \label{eq:isom1}
\end{equation*}
and, more generally,
\begin{equation}
(L/K,\sigma,a,v)\cong (L/K,\sigma,b,w) \iff a/b\in N_{L/K}(L).
  \label{eq:isom}
\end{equation}
\subsection{Generalized Cyclic Algebras}

Suppose $A$ contains a subfield $L$ that is cyclic over $K$ 
but not necessarily a maximal subfield.
Then $A$ is called a {\em generalized cyclic algebra}
(cf.\ generalized crossed products in Kursov-Yanchevski\u\i\ \cite{kursov-yanch-transl} 
or Tignol \cite{tignol:gen-cr-prod}).
Let $\sigma$ generate $\Gal(L/K)$, let $[L:K]=n_0$ and let $B$ denote the centralizer $Z_A(L)$ of $L$ in $A$.
By the Skolem-Noether Theorem there is an element $v\in A^*$ such that
$\Inn(v)|_L=\sigma$.
For any such $v$ we have $v^{n_0}\in B^*$ and
$$ A=\oplus_{i=0}^{n_0-1} Bv^i$$
as a $K$-space.
Since $B$ is the centralizer of $L$,
we have $vBv^{-1}=B$.
Setting $\wt\sigma:=\Inn(v)|_B$ and $a:=v^{n_0}$ we write
\begin{equation}
  A=(B/K,\wt\sigma,a,v).
  \label{eq:SN}
\end{equation}
Conversely, for any extension $\wt\sigma$ of $\sigma$ to $B$
there are $v\in A^*$ and $a\in B^*$ such that \eqref{eq:SN} holds.
If $\wt\sigma$ is fixed then 
\begin{equation}
(B/K,\wt\sigma,a,v)\cong (B/K,\wt\sigma,b,w) \iff a/b\in N_{L/K}(L)
  \label{eq:gen-isom}
\end{equation}
($a/b$ lies in $K$ if $a$ and $b$ arise in this way).

\subsection{Bicyclic Algebras}\label{sec:bicyclic}

We will use for bicyclic algebras the notation that was 
introduced in Amitsur-Saltman \cite[\S1]{amitsur-saltman:gen-abel-cr-prod}
for the more general ``abelian crossed products''
(see also Jacobson \cite[\S4.6, pp.174]{jacobson:fin-dim-div-alg}).
Suppose $A$ is bicyclic and let $F$ be a maximal subfield bicyclic over $K$.
Let $G:=\Gal(F/K)=G_1\times G_2$, $G_i=\gen{\sigma_i}$
and $|G_i|=n_i$.
We denote by $F_i$ the fixed field of 
$\sigma_i$
and by $N_i$ the norm map of the extension $F/F_i$.
Clearly, $F/F_i$ and $F_{3-i}/K$ each have group $G_i$.
By the Skolem-Noether Theorem there are elements $z_1,z_2\in A^*$ such that 
$\Inn(z_i)|_F=\sigma_i$ for $i=1,2$. 
For any such $z_1,z_2$ we have 
$$A=\bigoplus_{i=0}^{n_1-1}\bigoplus_{j=0}^{n_2-1} Fz_1^iz_2^j$$
as a $K$-space.
With the action of $z_i$ on $F$ being fixed,
the algebra $A$ is determined up to isomorphism by the elements
$$b_1:=z_1^{n_1}, \quad b_2:=z_2^{n_2}, \quad u:=z_2z_1z_2^{-1}z_1^{-1}.$$
We write
$$A=(F/K,z,u,b)$$
where $z=(z_1,z_2)$ and $b=(b_1,b_2)$.
The elements $b_1,b_2$ and $u\in F^*$ satisfy the relations
\begin{equation}
 N_1(u)=\frac{\sigma_2(b_1)}{b_1}, \quad N_2(u)=\frac{b_2}{\sigma_1(b_2)}, \quad b_i\in F_i
  \label{eq:nec-suf}
\end{equation}
(cf.\ \cite[Lemma 1.2]{amitsur-saltman:gen-abel-cr-prod}),
and \eqref{eq:nec-suf} are also sufficient for given elements $b_1,b_2,u\in F^*$
to define a bicyclic algebra
(cf.\ \cite[Theorem 1.3]{amitsur-saltman:gen-abel-cr-prod}).
If $u=1$ then \eqref{eq:nec-suf} imply $b_1,b_2\in K$ and we have
a canonical isomorphism
\begin{equation}
 (F/K,z,1,b)\cong (F_2/K,\sigma_1,b_1,z_1)\otimes(F_1/K,\sigma_2,b_2,z_2) 
  \label{eq:decomp}
\end{equation}
that identifies $F$ with $F_2\otimes F_1$ 
and the $z_i$ on both sides, respectively. 
In particular, we conclude that $(F/K,z,1,1)$ is split.

\begin{theorem}\label{thm:split}
$A=(F/K,z,u,b)$ is split if and only if there are $x_1,x_2\in F^*$ such that 
\begin{equation}
  N_1(x_1)=b_1,\quad N_2(x_2)=b_2, \quad \frac{\sigma_2(x_1)}{x_1}\frac{x_2}{\sigma_1(x_2)}=u.
  \label{eq:split}
\end{equation}
\end{theorem}
\begin{proof}
Since $(F/K,z,1,1)\cong M_{n}(K)$ the statement is a special case of \cite[Theorem 1.4]{amitsur-saltman:gen-abel-cr-prod}.
\end{proof}

\section{The Splitting Problem}\label{sec:sp}

Let $A$ be a central-simple $K$-algebra of degree $n$.

\begin{spp}
  Decide whether $A$ is split and, if so, compute a splitting of $A$,
  i.e.\ a $K$-isomorphism $A\ra M_n(K)$.
\end{spp}

For cyclic algebras the splitting problem quite obviously reduces to the solution of a norm equation.
The point of this section is to show the same for bicyclic algebras.
However, we start with the details of the cyclic case.

\begin{alg}[Splitting of cyclic algebra]\label{alg:split-case}Let
a cyclic algebra $A=(L/K,\sigma,a,v)$ of degree $n$ be given.
The splitting problem for $A$ is solved as follows.
\begin{enumerate}
\item Fix a $K$-embedding $\psi : L \ra M_n(K)$.
\item Compute $X\in M_n(K)$ such that 
  $\Inn(X)|_{\psi(L)}=\psi\sigma\psi^{-1}$
  and set $b:=X^n$. Then we have $b\in K$.
\item Solve the norm equation $N_{L/K}(x)=a/b$ for $x\in L$.
  If there is no solution then $A$ is not split, otherwise
  $\psi$ is extended to a splitting $A\ra M_n(K)$ by mapping $v$ to $\psi(x)X$.
\end{enumerate}
\end{alg}
\begin{proof}
Step 1 amounts to computing the minimal polynomial over $K$ of a primitive element of $L$.
The matrix $X$ in step 2 exists by the Skolem-Noether Theorem (cf.\ \S\ref{sec:cyclic})
and its computation is a linear problem.
Moreover, we have $b=X^n\in K$ and $M_n(K)=(\psi(L)/K,\psi\sigma\psi^{-1},b,X)$.
Step~3: By \eqref{eq:isom}, $A$ is split if and only if $a/b\in N_{L/K}(L)$.
If $x\in L$ is a solution to the equation $N_{L/K}(x)=a/b$ then 
$(\psi(x)X)^n=\psi(N_{L/K}(x))b=a$.
This shows that mapping $v$ to $\psi(x)X$ indeed defines an extension of $\psi$ to $A$. 
\end{proof}

\begin{remark}
  \label{rem:split-11}
  Algorithm \ref{alg:split-case} can be used to compute a splitting of the bicyclic algebra $(F/K,z,1,1)$.
Indeed, we take the canonical isomorphism
$$(F/K,z,1,1)\ra (F_2/K,\sigma_1,1,z_1)\otimes(F_1/K,\sigma_2,1,z_2)$$ 
from \eqref{eq:decomp},
compute splittings $(F_{3-i}/K,\sigma_i,1,z_i)\ra M_n(K)$ 
with Algorithm \ref{alg:split-case},
and finally compose with an isomorphism $M_n(K)\otimes M_n(K)\ra M_{n^2}(K)$.\end{remark}

Now let $A=(F/K,z,u,b)$ be an arbitrary bicyclic algebra. 
By Theorem \ref{thm:split}, $A$ is split if and only if 
the system of equations \eqref{eq:split} has a solution.  
In fact, any solution $(x_1,x_2)$ to \eqref{eq:split}
defines an isomorphism 
$$(F/K,z,u,b)\ra (F/K,w,1,1),\quad z_i\mapsto x_iw_i$$
(cf.\ the proof of \cite[Theorem 1.4]{amitsur-saltman:gen-abel-cr-prod}).
Together with Remark \ref{rem:split-11} we get a splitting of $A$. 
It remains to solve the system \eqref{eq:split},
and this covers the rest of the section.

\begin{lemma}[Bicyclic Hilbert 90]
  If $x_1,x_2\in F^*$ satisfy
  $$N_1(x_1)=1,\quad N_2(x_2)=1 \qt{and} \frac{\sigma_2(x_1)}{x_1}\frac{x_2}{\sigma_1(x_2)}=1$$
  then there is $y\in F^*$ with $x_i=\frac{\sigma_i(y)}{y}$ for $i=1,2$.
  \label{lem:bic-h90}
\end{lemma}
\begin{proof}
  In fact a stronger statement holds:
  if elements $x_1,\ldots,x_r\in F^*$ satisfy $N_i(x_i)=1$
  and $\frac{\sigma_i(x_j)}{x_j}=\frac{\sigma_j(x_i)}{x_i}$ for all $1\leq i,j\leq r$
  then $y\in F^*$ exists with $x_i=\frac{\sigma_i(y)}{y}$ for all $1\leq i\leq r$
  (cf.\ \cite[Proposition 4.6.30, p.179]{jacobson:fin-dim-div-alg}).
  However, we give a shorter proof for $r=2$.

  Let $x_1=\frac{\sigma_1(y_1)}{y_1}$ with $y_1\in F^*$, 
  by Hilbert's Theorem 90.
  Then
  $$\frac{\sigma_1(x_2)}{x_2}=\frac{\sigma_2(x_1)}{x_1}=\frac{\sigma_2\sigma_1(y_1)\cdot y_1}{\sigma_2(y_1)\cdot\sigma_1(y_1)}
  =\frac{\sigma_1(\frac{\sigma_2(y_1)}{y_1})}{\frac{\sigma_2(y_1)}{y_1}},$$
  hence $x_2=c\cdot\frac{\sigma_2(y_1)}{y_1}$ for some $c\in F_1^*$.
  It follows $N_2(c)=1$. 
  Let $c=\frac{\sigma_2(y_2)}{y_2}$ with $y_2\in F_1^*$, by Hilbert's Theorem 90.
  Defining $y:=y_1y_2$ we get $x_i=\frac{\sigma_i(y)}{y}$ for $i=1,2$.
\end{proof}

\begin{prop}\label{prop:S}
  Suppose \eqref{eq:split} has a solution $(x_1,x_2)$.
  Then the set of all solutions is
  $$S:=\sett{(x_1\frac{\sigma_1(y)}{y},x_2\frac{\sigma_2(y)}{y})}{y\in F^*}.$$
  In particular, for any $x_1'\in F^*$ with $N_1(x_1')=b_1$ there is $x_2'\in F^*$ such that $(x_1',x_2')\in S$.
\end{prop}
\begin{proof}
  An easy calculation shows that any $(x'_1,x'_2)\in S$ solves \eqref{eq:split}.
  For the converse apply Lemma \ref{lem:bic-h90}.
  The second statement is another application of Hilbert's Theorem~90.
\end{proof}

\begin{alg}[Solution to \eqref{eq:split}]
  \label{alg:split}
If a solution to \eqref{eq:split} exists
then one is found performing the following steps.
  Conversely, if all steps have a solution then
  the resulting $(x_1,x_2)$ is a solution to \eqref{eq:split}.
\begin{enumerate}
\item Solve $N_1(x_1)=b_1$ for $x_1\in F^*$.
\item Solve 
  $$\frac{\sigma_1(x'_2)}{x'_2} = u\frac{x_1}{\sigma_2(x_1)}
\qt{for $x'_2\in F^*$.}$$
\item Solve $N_2(x''_2)=b_2 N_2(x'_2)$ for $x''_2\in F_1^*$.
\item Define $x_2:={x'_2}^{-1} x''_2$.
\end{enumerate}
\end{alg}
\begin{proof}
  A straight-forward calculation verifies that any $(x_1,x_2)$ computed by these steps
  is a solution to \eqref{eq:split}.
  Conversely, suppose \eqref{eq:split} has a solution and show each step has a solution.
  Step 1 is obvious.
  Step 2 has a solution by Hilbert's Theorem 90 because, using \eqref{eq:nec-suf},
  $$N_1(u\frac{x_1}{\sigma_2(x_1)})=N_1(u)\frac{b_1}{\sigma_2(b_1)}=1.$$
  Step 3: Since \eqref{eq:split} has a solution,  
  Proposition \ref{prop:S} shows the existence of an element $x_2\in F^*$ with $(x_1,x_2)\in S$.
  Setting $x_2'':=x_2x'_2$, \eqref{eq:split} implies $\sigma_1(x''_2)/{x''_2}=1$ and $N_2(x''_2)=b_2N_2(x'_2)$.
  This shows that step 3 has a solution.
\end{proof}

  Algorithm \ref{alg:split} reduces the splitting of a bicyclic algebra
  to two consecutive ({\em not simultaneous}) norm equations;
  the rest (step 2) is linear algebra.
  Note that the first norm equation (step 1) lives in the larger fields $F/F_1$
  whereas the second one (step 3) lives in $F_1/F$.
    Algorithms for norm equations over number fields can be found in
    Simon \cite{simon:neq} and the references cited therein.

  \begin{rem}
    The splitting problem has two parts: first, to decide whether the algebra
    is split, and second, to compute a splitting.  One might be tempted to think
    that it suffices for the first part to decide the solvability of
    norm equations, and that solutions are required only for the
    second part.  Indeed, this is true for cyclic algebras because only one
    norm equation appears.  For bicyclic algebras, however, two norm
    equations occur in Algorithm \ref{alg:split} 
    and the second one is built from a solution to the first.  
    Thus, at least the first norm equation has to be solved explicitly
    for any part of the splitting problem.
  \end{rem}

  \comment{
  \begin{rem}
    Since the known algorithms involve factoring algebraic integers.
    Norm equations are thus considerably ``hard'' to solve.

  \end{rem}
  }

\section{The Isomorphism Problem}\label{sec:ip}

Let $A_1$ and $A_2$ be central-simple $K$-algebras of degree $n$.

\begin{ip}
  Decide whether $A_1$ and $A_2$ are $K$-isomorphic and, if so,
  compute a $K$-isomorphism between them.
\end{ip}

We show for general $A_1$ and $A_2$ 
how the isomorphism problem reduces to the splitting problem
(see Algorithm \ref{alg:A-Aop} below).
This is due to the equivalence \eqref{eq:op}. 
When specializing thereafter to cyclic algebras,  
the isomorphism problem eventually reduces to norm equations.

\begin{remark}\label{rem:emb}
A $K$-algebra isomorphism $\varphi:A_1\otimes A_2^\circ\ra C$ is equivalent
to a pair $(\varphi_1,\varphi_2)$ where
$\varphi_1:A_1\ra C$ is a $K$-embedding
and 
$\varphi_2:A_2\ra C$ is a $K$-anti-embedding
such that $\varphi_1(A_1)$ is the centralizer of $\varphi_2(A_2)$ in $C$.
Of course, one obtains $\varphi_1,\varphi_2$ from $\varphi$ by composing $\varphi$ 
with the canonical embedding $\eps_1: A_1\ra A_1\otimes A_2^\circ$ and canonical anti-embedding $\eps_2: A_2\ra A_1\otimes A_2^\circ$, respectively.
\end{remark}

\begin{alg}
  \label{alg:A-Aop}
  Suppose a $K$-isomorphism $\varphi:A_1\otimes A_2^\circ\ra M_{n^2}(K)$ is 
given in the form of a pair $(\varphi_1,\varphi_2)$ as in Remark~\ref{rem:emb}.
Then $K$-isomorphisms $\chi: A_1\ra A_2$ and $\chi': A_2\ra A_1$ are computed as follows. 
\begin{enumerate}
\item Fix a $K$-basis of $A_1$ and with respect to this basis
  identify $M_{n^2}(K)$ with $\End_K(A_1)$.
\item Compute $X\in M_{n^2}(K)$ such that $\Inn(X)\circ\varphi_1$ is
the left-regular representation of $A_1$.
\item Set $\varphi'_2 := \Inn(X)\circ\varphi_2$.
  Then $\varphi'_2(A_2)=\rho(A_1)$, where $\rho$ is the right-regular representation of~$A_1$.
\item Define $\chi:={\varphi'_2}^{-1}\circ\rho$ and $\chi':=\rho^{-1}\circ{\varphi'_2}$.
\end{enumerate}
\end{alg}
\begin{proof}
The left-regular representation
$$\lambda : A_1\ra\End_K(A_1), \quad a\mapsto \lambda_a,\quad \lambda_a(x) := ax$$ 
is a $K$-algebra embedding 
and the right-regular representation  
$$\rho : A_1\ra\End_K(A_1), \quad a\mapsto \rho_a,\quad \rho_a(x) := xa$$ 
is a $K$-algebra anti-embedding.
The matrix $X$ in step 2 exists by the Skolem-Noether Theorem and
its computation is a linear problem.
We get
\begin{multline*}
 \varphi'_2(A_2)=
Z_{M_{n^2}(K)}(X\varphi_1(A_1)X^{-1})
=Z_{M_{n^2}(K)}(\lambda(A_1))\\
=\rho(A_1). 
\end{multline*}
Since $\rho$ and $\varphi'_2$ are both anti,
$\chi$ and $\chi'$ as defined in step 4 are isomorphisms.
\end{proof}

\begin{rem}
The reduction of the isomorphism problem for $A_1$ and $A_2$
to the splitting problem for $A_1\otimes A_2^\circ$
as in Algorithm \ref{alg:A-Aop} 
is accompanied by an increase in the algebra degree from $n$ to $n^2$. 
\end{rem}

Now we turn to the case when $A_1$ and $A_2$ are both cyclic, say
$$A_1=(L_1/K,\sigma_1,a_1,v_1)\qt{and} A_2=(L_2/K,\sigma_2,a_2,v_2).$$
We distinguish two special cases and the general case.

\begin{spcase}
$L_1\cong_K L_2$. 
Let $\chi : L_1\ra L_2$ be a $K$-isomorphism.
Obviously, $\chi\sigma_1\chi^{-1}=\sigma_2^i$ for some $i$ relatively prime to $n$. 
Replacing $v_2$ with $v_2^i$
we can assume by \eqref{eq:k} that $\chi\sigma_1\chi^{-1}=\sigma_2$.
Then, by \eqref{eq:isom},
$$A_1\cong A_2 \iff a_1/a_2\in N_{L_2/K}(L_2).$$
This equivalence is ``constructive'':
if $x\in L_2$ is an element with $N_{L_2/K}(x)=a_1/a_2$ 
then mapping $v_1\mapsto xv_2$ extends $\chi$ to a $K$-isomorphism $A_1\ra A_2$.
The isomorphism problem for $A_1$ and $A_2$ is thus reduced to
finding a solution to a norm equation in the field extension $L_2/K$.
\end{spcase}

\begin{spcase}
$L_1/K$ and $L_2/K$ are linearly disjoint, i.e.\
$F:=L_1\otimes L_2$ is a field.
By \eqref{eq:decomp}, $A_1\otimes A_2^\circ$ is isomorphic to the bicyclic algebra
$C:=(F/K,z,1,b)$
with $z=(z_1,z_2)$ and $b=(a_1,a_2^{-1})$.
(Here, we regard $\sigma_1,\sigma_2$ as automorphisms of $F$.)
Indeed, the canonical isomorphism $\varphi:A_1\otimes A_2^\circ\ra C$ maps
$\lambda v_1^i\otimes 1\mapsto \lambda z_1^i$ for all $\lambda\in L_1$
and $1\otimes\lambda v_2^i\mapsto z_2^{-i}\lambda$ for all $\lambda\in L_2$.
Thus, Algorithm \ref{alg:A-Aop} reduces the isomorphism problem for $A_1$ and $A_2$ 
to the splitting problem for $C$.
Since $C$ is bicyclic, this further reduces to norm equations by \S\ref{sec:sp}.
\end{spcase}

\begin{rem}
If the degree $n$ is prime then we are in one of the special cases.
\end{rem}

\begin{gcase}
Regarding $L_1,L_2$ as subfields of some common overfield,
we set $L_0:=L_1\cap L_2$ and consider the centralizers $B_1:=Z_{A_1}(L_0)$ and $B_2:=Z_{A_2}(L_0)$.
If $A_1\cong A_2$ there is by the Skolem-Noether Theorem a $K$-isomorphism
$\chi:A_1\ra A_2$ with $\chi(L_0)=L_0$ and $\chi|_{L_0}=\id_{L_0}$.
For any such $\chi$ we have $\chi(B_1)=B_2$,
i.e.\ the restriction $\chi|_{B_1}$ 
is an $L_0$-isomorphism $B_1\ra B_2$. 

We therefore solve the general case by starting with
the isomorphism problem for $B_1$ and $B_2$ over $L_0$
(which is special case 2).
If it has no solution then $A_1\not\cong A_2$. 
Assume otherwise and suppose an $L_0$-isomorphism $\chi_0:B_1\ra B_2$ is computed.
We identify $B_1$ and $B_2$ under $\chi_0$ and simply write $B$ for both of them.
By \S\ref{sec:bicyclic},
there are elements $w_1\in A_1^*, w_2\in A_2^*$ and $b_1,b_2\in B^*$ such that
$$A_1=(B/K,\wt\sigma,b_1,w_1)\qt{and} A_2=(B/K,\wt\sigma,b_2,w_2)$$
with the same $\wt\sigma$ for both algebras.
Since $w_1,w_2$ are predicted by the Skolem-Noether Theorem,
their computation is a linear problem.
We proceed as in special case 1 but for generalized cyclic algebras.
By \eqref{eq:gen-isom},
$$A_1\cong A_2 \iff b_1/b_2\in N_{L_0/K}(L_0).$$
This equivalence is ``constructive'':
if $x\in L_0$ is an element with $N_{L_0/K}(x)=b_1/b_2$ 
then mapping $w_1\mapsto xw_2$ defines a $K$-isomorphism $A_1\ra A_2$.
Thus, also the general case is reduced to norm equations.
\end{gcase}

\section{Extending Field Automorphisms to Simple Algebras}\label{sec:ext}

Let $A$ be a central-simple $K$-algebra 
and let $\sigma$ be an automorphism of $K$ of finite order.
\begin{ep}
Decide whether $\sigma$ extends to an automorphism of $A$ and, if so, compute an extension.
\end{ep}
It is convenient to reformulate this problem using the algebra $\pre{\sigma^{-1}}{A}$ which is obtained from $A$ by redefining the $K$-action as
$\lambda\circ a:=\sigma^{-1}(\lambda)a$ for all $\lambda\in K$.
Then $\sigma$ extends to $A$ if and only if $A$ and $\pre{\sigma^{-1}}{A}$ are isomorphic as $K$-algebras.
In fact, any $K$-algebra isomorphism $\pre{\sigma^{-1}}{A}\ra A$ becomes an extension of $\sigma$
after identifying $\pre{\sigma^{-1}}{A}$ as a ring with $A$.
The extension problem is therefore just a special case of the isomorphism problem.

If $A$ is cyclic then $\pre{\sigma^{-1}}{A}$ is also cyclic,
hence, by the results of the preceding sections,
the extension problem for cyclic algebras reduces to norm equations.

\begin{remark}\label{rem:ring-id}
Let $A=(L/K,\tau,a,v)$ and identify $\pre{\sigma^{-1}}{A}$ as a ring with $A$.
For any extension of $\sigma$ to $L$ (call it also $\sigma$),
we have
$$\pre{\sigma^{-1}}{A}=(\sigma L/K,\sigma\tau\sigma^{-1},\sigma a,v).$$
The ring identity map $A\ra \pre{\sigma^{-1}}{A}$ is defined by $\sigma: L\ra \sigma L$ and $v\mapsto v$.
\end{remark}

We finish with a detailed example for the solution of the extension problem.
Let $K$ be the cubic number field 
\begin{equation*}
K=\Q(\alpha), \quad \Irr(\alpha,\Q)=x^3+x^2-2x-1,
\end{equation*}
of discriminant $49$
(the maximal real subfield of the $7$-th cyclotomic field).
We have $\Gal(K/\Q)=\gen{\sigma}$ with 
$$\sigma(\alpha)= -\alpha^2-\alpha+1.$$
Let $L$ be the cubic extension
\begin{equation*}
L=K(\theta), \quad \Irr(\theta,K)=x^3+(\alpha-2)x^2+(-\alpha-1)x+1,
\end{equation*}
which is cyclic and has $\Gal(L/K)=\gen{\tau}$ with
\begin{equation*}
  \tau(\theta) = -\theta^2 + (-\alpha + 1)\theta + 2.
\end{equation*}
We will solve the extension problem for $\sigma$ and the cyclic division algebra
$$D=(L/K,\tau,a,v), \quad a=2(\alpha^2-\alpha-2).$$
In order to do so 
we solve the isomorphism problem for $D$ and $\pre{\sigma^{-1}}D$.
According to Remark \ref{rem:ring-id}
we have $$\pre{\sigma^{-1}}D=(\sigma L/K,\sigma\tau\sigma^{-1},\sigma a,v),$$
where $\sigma L=K(\eta)$ with 
$$\Irr(\eta,K)=x^3+(-\alpha^2-\alpha-1)x^2+(\alpha^2+\alpha-2)x+1,$$
$\sigma\tau\sigma^{-1}(\eta) = -\eta^2 + (\alpha^2 + \alpha)\eta + 2$
and $\sigma(a) =2(\alpha^2+2\alpha-1)$. 
Of course, $\Irr(\eta,K)$ is obtained from $\Irr(\theta,K)$ by applying $\sigma$ to each coefficient.

Since $L_1/K$ and $L_2/K$ are linearly disjoint, we proceed as in special case 2
and consider the bicyclic algebra 
$$ C:=(F/K,z,1,b)$$
with $F=L\otimes\sigma L, \sigma_1=\tau\otimes\id, \sigma_2=\id\otimes\sigma\tau\sigma^{-1}, F_1=\sigma L, F_2=L, z=(z_1,z_2)$
and $b=(a,\sigma a^{-1})$.
The canonical isomorphism $D\otimes \pre{\sigma^{-1}}D^\circ\ra C$
maps 
$\lambda v^i\otimes 1\mapsto \lambda z_1^i$
and 
$1\otimes\sigma(\lambda) v^i\mapsto z_2^{-i}\sigma(\lambda)$
for all $\lambda\in L$.
We compute a splitting of $C$ following~\S\ref{sec:sp}.

First apply Algorithm \ref{alg:split}.
Step 1. 
The computer algebra system MAGMA \cite{magma} gives 
\tiny
\begin{equation*}
  \begin{split}
    x_1 = &\frac{1}{2}\left((-7\alpha^2 + 9\alpha + 4)+2(-2\alpha^2 + 6\alpha - 1)\eta+(\alpha^2 - 6\alpha + 4)\eta^2\right.\\
    &+(14\alpha^2 - 12\alpha - 9)\theta+(20\alpha^2 - 7\alpha - 14)\eta\theta+(-12\alpha^2 + 3\alpha + 12)\eta^2\theta\\
    &\left.+(-3\alpha^2 + 3\alpha + 2)\theta^2+(-7\alpha^2 - \alpha + 5)\eta\theta^2+(4\alpha^2 - 5)\eta^2\theta^2\right). 
  \end{split}
  \label{eq:sol1}
\end{equation*}
\normalsize
Step 2. As a solution to a linear equation system one finds
\tiny
\begin{equation*}
  \begin{split}
    x'_2 = &\frac{1}{28}
    \left( 
    (81\alpha^2 - 261\alpha - 173)
    + (194\alpha^2 + 428\alpha + 132)\eta 
    + (-70\alpha^2 - 98\alpha - 21)\eta^2 
    \right.\\
    &
    +(- 19\alpha^2 - 120\alpha + 8)\theta
    + (287\alpha^2 - 266\alpha - 196)\eta\theta 
    + (-92\alpha^2 + 130\alpha + 45)\eta^2\theta 
    \\
    &\left. 
    + (7\alpha^2 + 70\alpha + 14)\theta^2
    + (-202\alpha^2 - 179\alpha - 24)\eta\theta^2 
    + (63\alpha^2 + 28\alpha + 7)\eta^2\theta^2 
    \right).
  \end{split}
\end{equation*}
\normalsize
\\
Step 3. Exceptionally in this example, the element
$$b_2N_2(x'_2) = \frac{1}{56}(-1601\alpha^2 +693\alpha+609)$$
lies in $K$. We compute as a cubic root~:
\begin{equation*}
  x''_2 = \frac{1}{14}(-19\alpha^2 - \alpha - 6).
\end{equation*}
\\
Step 4. Finally, we get
\tiny
\begin{equation*}
  \begin{split}
    x_2 = &\frac{1}{14}
    \left( 
    (-8\alpha^2 - 3\alpha + 10)
    + (-2\alpha^2 - 13\alpha - 1)\eta 
    + (2\alpha^2 + 6\alpha - 6)\eta^2 
    \right. \\ 
    &
    + (-3\alpha^2 + 5\alpha + 2)\theta 
    + (6\alpha^2 - 10\alpha + 10)\eta\theta 
    + ( -2\alpha^2 + \alpha - 1)\eta^2\theta 
    \\
    &\left.
    + (2\alpha^2 - \alpha - 6)\theta^2
    + (\alpha^2 + 10\alpha - 3)\eta\theta^2 
    + (-\alpha^2 - 3\alpha + 3)\eta^2\theta^2 
    \right).
  \end{split}
\end{equation*}
\normalsize
This completes Algorithm \ref{alg:split}.
On a $32$ bit machine with an 800~MHz processor
the computation time for the norm equation in step~1 was less than ten minutes
and was negligible for steps 2---4.

At this point we have shown that $C$ is split,
hence $\sigma$ extends to $D$.
We continue to compute an extension of $\sigma$.
Since $(x_1,x_2)$ is a solution to \eqref{eq:split},
we obtain an isomorphism $\varphi : C\ra M_{9}(K)$
as it is described after Algorithm \ref{alg:split-case}.
  Using Algorithm \ref{alg:A-Aop} 
  we get the $K$-isomorphism $\chi':\pre{\sigma^{-1}}D\ra D$
  defined by
  \tiny
  \begin{multline*}
    \chi'(\eta) := \frac{1}{673}\cdot(1,\theta,\theta^2)\cdot\\
    \chieta\cdot\left(\begin{matrix}1\\v\\v^2\end{matrix}\right). 
  \end{multline*}
  \normalsize and \tiny
  \begin{equation*}
    \chi'(w) := (1,\theta,\theta^2)\cdot\chiw\cdot\left(\begin{matrix}1\\v\\v^2\end{matrix}\right)
      = \left((\alpha^2+\alpha)+(-\alpha+1)\theta-\theta^2\right)v.
  \end{equation*}
  \normalsize
The resulting extension $\wt\sigma$ of $\sigma$ is 
$$\wt\sigma: D\ra D, \quad \theta\mapsto\chi'(\eta), \quad v\mapsto\chi'(w).$$

\begin{rem}
The example is taken from the paper \cite{hanke:laurent-noncr}
which uses $D$ and $\wt\sigma$ to construct a noncrossed product division algebra
in the form of the {\em twisted Laurent series ring} 
$D(\!(\x;\wt\sigma)\!)$,
i.e.\ the ring of all formal series $\sum_{i\geq k} d_i\x^i, k\in\Z,$ with multiplication of monomials
$d\x^i\cdot d'\x^j=d\wt\sigma^i(d')\x^{i+j}$.
The ring $D(\!(\x;\wt\sigma)\!)$ is a division algebra of degree $9$ over the power series field $\Q(\!(\bf t)\!)$.
As shown in \cite{hanke:laurent-noncr},
$D$ does not contain a maximal subfield that is Galois over $\Q$,
and this property implies that $D(\!(\x;\wt\sigma)\!)$ is a noncrossed product.
The above computation explains how the $\wt\sigma$ given in \cite{hanke:laurent-noncr} was found.
(Note that $D$ in \cite{hanke:laurent-noncr}
is accidentally defined with $-a$ instead of $a$.
This sign error is corrected in the version at arXiv:math/0703038.)
\end{rem}

\section{Acknowledgements}
I would like to thank the department of mathematics at 
Universidad Nacional Aut\'onoma de M\'exico, M\'exico, D.F.,
and in particular my host Jos\'e Antonio de la Pe\~na
for their kind invitation and support.
I am also indebted to Adrian Wadsworth and the department of mathematics
at the University of California at San Diego.
It was at San Diego where the first results of this work were obtained.
The financial support for this visit
from the German Academic Exchange Service (DAAD) under grant D/02/00701
is acknowledged.
Special thanks are due to the MAGMA group for making available the system to me
and for the advice on its use.
The helpful comments of the referees were gratefully appreciated.
\bibliographystyle{abbrv}

\begin{thebibliography}{10}

\bibitem{amitsur-saltman:gen-abel-cr-prod}
S.~Amitsur and D.~Saltman.
\newblock Generic abelian crossed products and $p$-algebras.
\newblock {\em J. Algebra}, 51:76--87, 1978.

\bibitem{magma}
W.~Bosma, J.~Cannon, and C.~Playoust.
\newblock The magma algebra system {I} : The user language.
\newblock {\em J. Symb. Comp.}, 24(3/4):235--265, 1997.
\newblock (see the Magma home page at
  http://www.maths.usyd.edu.au:8000/u/magma/)

\bibitem{kash}
M.~Daberkow, C.~Fieker, J.~Kl{\"u}ners, M.~Pohst, K.~Roegner, and K.~Wildanger.
\newblock {KANT V4}.
\newblock {\em J. Symbolic Comp.}, 24:267--283, 1997.

\bibitem{graaf-schicho}
W.~de~Graaf, M.~Harrison, J.~P{\'\i}lnikov\'a, and J.~Schicho.
\newblock {A Lie algebra method for rational parametrization of Severi-Brauer
  surfaces.}
\newblock {\em J. Algebra}, 303(2):514--529, 2006.

\bibitem{eberly:decomp-nf}
W.~Eberly.
\newblock {Decomposition of algebras over finite fields and number fields.}
\newblock {\em Comput. Complexity}, 1(2):183--210, 1991.

\bibitem{haile:useful}
D.~Haile.
\newblock A useful proposition for division algebras of small degree.
\newblock {\em Proc. Amer. Math. Soc.}, 106(2):317--319, 1989.

\bibitem{hanke:laurent-noncr}
T.~Hanke.
\newblock A twisted {L}aurent series ring that is a noncrossed product.
\newblock {\em Israel J. Math.}, 150:199--204, 2005.

\bibitem{jacobson:fin-dim-div-alg}
N.~Jacobson.
\newblock {\em Finite Dimensional Division Algebras over Fields}.
\newblock Springer-Verlag, Berlin, 1996.

\bibitem{kursov-yanch-transl}
V.~V. Kursov and V.~I. Yanchevski\v\i.
\newblock Crossed products of simple algebras and their automorphism groups.
\newblock {\em Amer. Math. Soc. Transl.}, 154(2), 1992.

\bibitem{pierce:ass-alg}
R.~Pierce.
\newblock {\em {Associative Algebras}}.
\newblock Springer-Verlag, New York, 1982.

\bibitem{reiner:max-orders}
I.~Reiner.
\newblock {\em {Maximal Orders}}.
\newblock Academic Press, London, 1975.

\bibitem{simon:neq}
D.~Simon.
\newblock Solving norm equations in relative number fields using {$S$}-units.
\newblock {\em Math. Comp.}, 71(239):1287--1305, 2002.

\bibitem{tignol:gen-cr-prod}
J.-P. Tignol.
\newblock Generalized crossed products.
\newblock S\'eminaire Math\'ematique (nouvelle s\'erie), No. 106, Universit\'e
  Catholique de Louvain, Louvain-la-Neuve, Belgium, 1987.

\end{thebibliography}
\def\cprime{$'$}

\end{document}